\documentclass[1p,12pt]{elsarticle}

\usepackage{graphicx}
\usepackage{amsmath,amsxtra,amssymb,latexsym, amscd,amsthm}

\theoremstyle{plain}
\newtheorem{thm}{Theorem}[section]

\theoremstyle{definition}

\theoremstyle{remark}

\numberwithin{equation}{section}
\numberwithin{figure}{section}
\numberwithin{table}{section}

\newcommand{\M}{\operatorname{M}}

\begin{document}

\begin{frontmatter}

\title{\textbf{A Generalization of Aztec Diamond Theorem, Part II}}

\author{Tri Lai\corref{cor1}\fnref{myfootnote1}}
\address{Institute for Mathematics and its Applications\\ University of Minnesota, Minneapolis, MN 55455}
\fntext[myfootnote1]{This research was supported in part by the Institute for Mathematics and its Applications with funds provided by the National Science Foundation (grant no. DMS-0931945).}
\cortext[cor1]{Corresponding author, email: tmlai@ima.umn.edu, tel: 612-626-8319}

\begin{abstract}
The author gave a proof of a generalization of the Aztec diamond theorem  for a family of $4$-vertex regions on the square lattice with southwest-to-northeast diagonals drawn in (\textit{Electron. J.  Combin.}, 2014) by using a bijection between tilings and non-intersecting lattice paths. In this paper, we use Kuo graphical condensation to give a new proof.
\end{abstract}


\end{frontmatter}

\section{Introduction}

A \emph{lattice} divides the plane into non-overlapped parts called \textit{fundamental regions}. A \textit{region} considered in this paper is a finite connected union of fundamental regions. We define a \textit{tile} to be the union of  any two fundamental regions sharing an edge. We are interested in how many different ways to cover a region by tiles such that there are no gaps or overlaps; such coverings are called \textit{tilings}.  We use the notation $\M(R)$ for the number of tilings of a region $R$.

\begin{figure}\centering
\includegraphics[width=7cm]{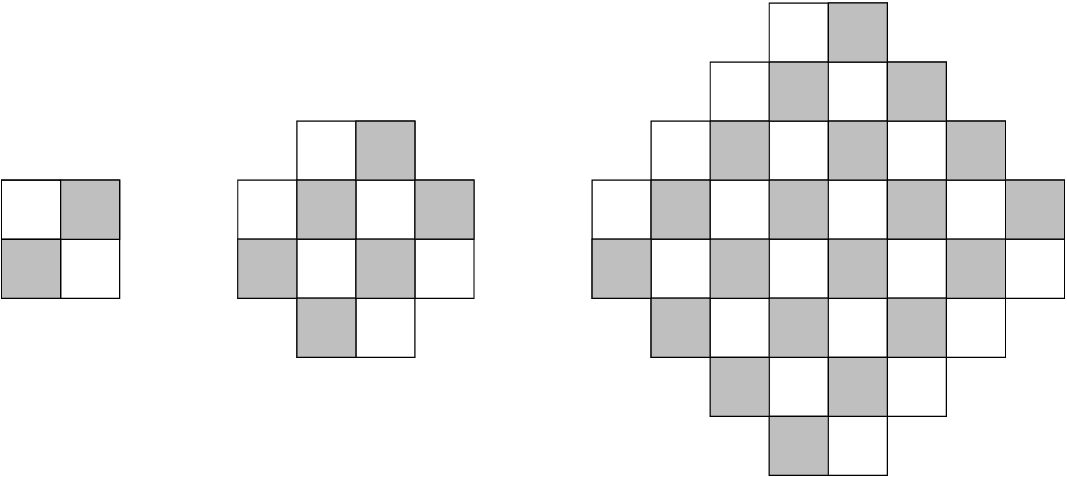}
\caption{From left to right, the Aztec diamonds of order $1$, $2$ and $4$.}
\label{aztecdiamond}
\end{figure}

The \textit{Aztec diamond} of order $n$ is the union of all the unit squares inside the contour $|x|+|y|=n+1$. Figure \ref{aztecdiamond} illustrates the Aztec diamonds of order $1,$ $2,$ and $4$ with a checkboard coloring.
The well-known Aztec diamond theorem by Elkies, Kuperberg, Larsen and Propp \cite{Elkies, Elkies2} states that the number of (domino) tilings of the Aztec diamond of order $n$ is equal to $2^{n(n+1)/2}$. 
The first four proofs of the theorem were presented in \cite{Elkies, Elkies2}, and  many further proofs followed (see e.g. \cite{Bosio, Brualdi, Eu, Fend, Kamioka, Kuo, propp}).

\begin{figure}
\centering
\includegraphics[width=10cm]{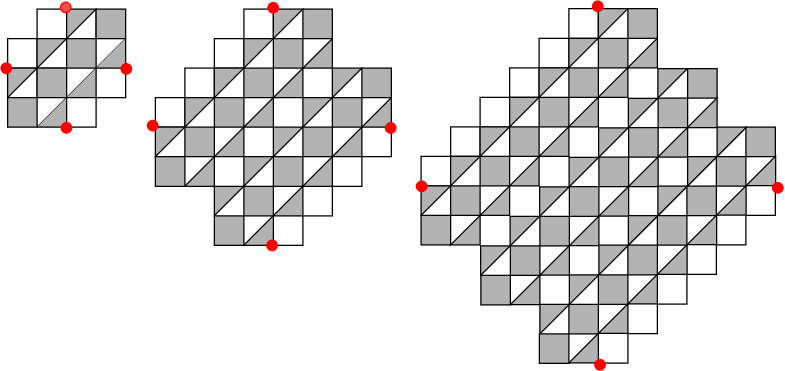}
\caption{From left to right, the Douglas regions of order $1$, $2$ and $3$.}
\label{douglas}
\end{figure}

Chris Douglas \cite{Doug} considered a variant of the Aztec diamond on the square lattice with every second southwest-to-northeast diagonal drawn in. The first three Douglas regions $D(n)$'s  are illustrated in Figure \ref{douglas}. More precisely, the four vertices of $D(n)$ (indicated by the dots in Figure \ref{douglas}) are always the vertices of a diamond with side-length $2n\sqrt{2}$. The northwest and the southeast boundaries of $D(n)$ are the same as that of the Aztec diamond of order $2n$, and the northeast and the southwest boundaries are two zigzag paths with steps of length 2. Douglas \cite{Doug} proved a conjecture posed by Propp that the region $D(n)$ has $2^{2n(n+1)}$ tilings.

We consider a family of regions first introduced in \cite{Tri}, which can be considered as a common generalization of the Aztec diamonds and Douglas regions.

From now on, the term ``diagonal" will be used to mean ``southwest-to-northeast lattice diagonal". 
The distances between any two successive drawn-in diagonals of the Douglas region $D(n)$ are all $\sqrt{2}$. Next, we consider the general situation when the distances  between two successive drawn-in diagonals are arbitrary.

\begin{figure}\centering
\setlength{\unitlength}{3947sp}%
\begingroup\makeatletter\ifx\SetFigFont\undefined%
\gdef\SetFigFont#1#2#3#4#5{%
  \reset@font\fontsize{#1}{#2pt}%
  \fontfamily{#3}\fontseries{#4}\fontshape{#5}%
  \selectfont}%
\fi\endgroup%
\resizebox{10cm}{!}{
\begin{picture}(0,0)%
\includegraphics{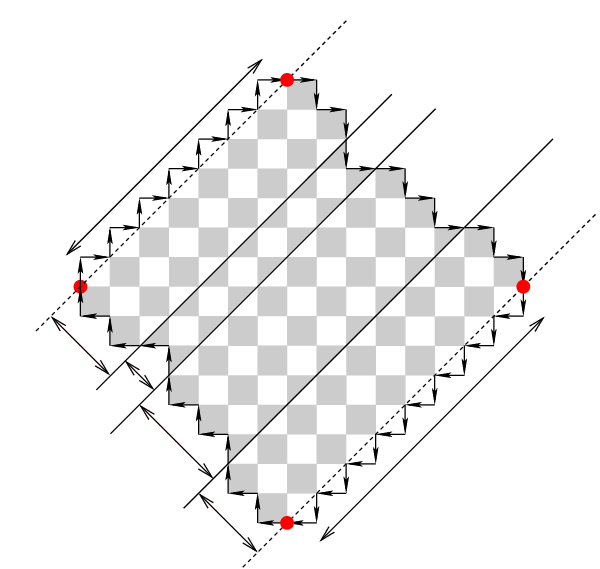}%
\end{picture}%

\begin{picture}(4794,4535)(924,-4323)
\put(5446,-1591){\makebox(0,0)[lb]{\smash{{\SetFigFont{12}{14.4}{\rmdefault}{\mddefault}{\itdefault}{$\ell'$}%
}}}}
\put(1800,-4216){\makebox(0,0)[lb]{\smash{{\SetFigFont{12}{14.4}{\rmdefault}{\mddefault}{\itdefault}{$\textbf{d}_4\frac{\sqrt{2}}{2}=\textbf{4}\frac{\sqrt{2}}{2}$}%
}}}}
\put(3136,-196){\makebox(0,0)[lb]{\smash{{\SetFigFont{12}{14.4}{\rmdefault}{\mddefault}{\itdefault}{$A$}%
}}}}
\put(3496, 36){\makebox(0,0)[lb]{\smash{{\SetFigFont{12}{14.4}{\rmdefault}{\mddefault}{\itdefault}{$\ell$}%
}}}}
\put(5347,-2182){\makebox(0,0)[lb]{\smash{{\SetFigFont{12}{14.4}{\rmdefault}{\mddefault}{\itdefault}{$B$}%
}}}}
\put(3221,-4308){\makebox(0,0)[lb]{\smash{{\SetFigFont{12}{14.4}{\rmdefault}{\mddefault}{\itdefault}{$C$}%
}}}}
\put(1704,-939){\makebox(0,0)[lb]{\smash{{\SetFigFont{12}{14.4}{\rmdefault}{\mddefault}{\itdefault}{$a=7$}%
}}}}
\put(4479,-3474){\makebox(0,0)[lb]{\smash{{\SetFigFont{12}{14.4}{\rmdefault}{\mddefault}{\itdefault}{$w=8$}%
}}}}
\put(600,-2821){\makebox(0,0)[lb]{\smash{{\SetFigFont{12}{14.4}{\rmdefault}{\mddefault}{\itdefault}{$\textbf{d}_1\frac{\sqrt{2}}{2}=\textbf{4}\frac{\sqrt{2}}{2}$}%
}}}}
\put(850,-3234){\makebox(0,0)[lb]{\smash{{\SetFigFont{12}{14.4}{\rmdefault}{\mddefault}{\itdefault}{$\textbf{d}_2\frac{\sqrt{2}}{2}=\textbf{2}\frac{\sqrt{2}}{2}$}%
}}}}
\put(1450,-3646){\makebox(0,0)[lb]{\smash{{\SetFigFont{12}{14.4}{\rmdefault}{\mddefault}{\itdefault}{$\textbf{d}_3\frac{\sqrt{2}}{2}=\textbf{5}\frac{\sqrt{2}}{2}$}%
}}}}
\put(1201,-2124){\makebox(0,0)[lb]{\smash{{\SetFigFont{12}{14.4}{\rmdefault}{\mddefault}{\itdefault}{$D$}%
}}}}
\end{picture}}
\caption{The region $D_{7}(4,2,5,4)$.}
\label{DouglasSchroder}
\end{figure}
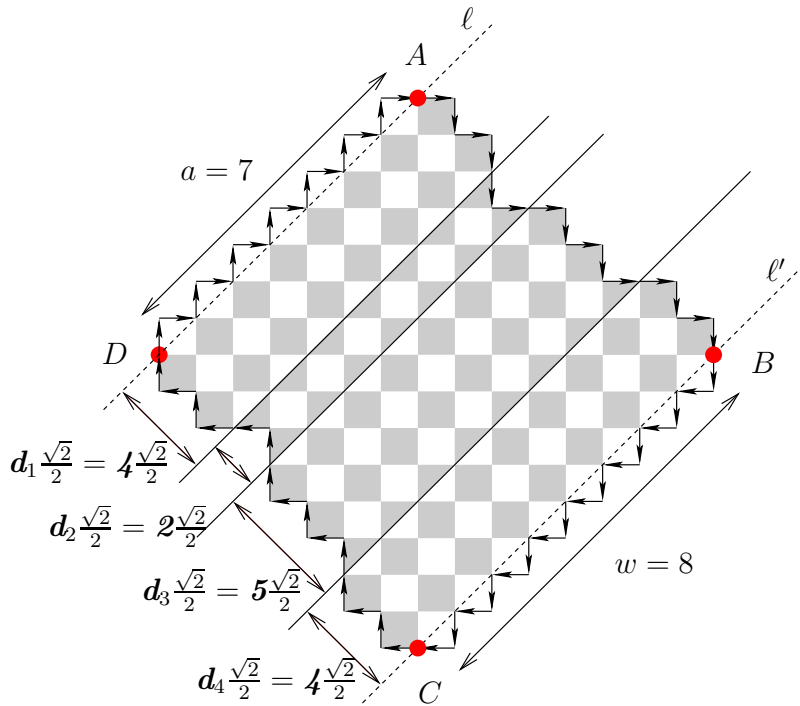

Suppose we have two lattice diagonals $\ell$ and $\ell'$ that are \textit{not} drawn-in diagonals, so that $\ell'$ is below $\ell$. Assume that $k-1$ diagonals have been drawn between $\ell$ and $\ell'$, with the distances between successive ones, starting from top, being $d_1\frac{\sqrt{2}}{2},d_2\frac{\sqrt{2}}{2},\dotsc,d_{k-1}\frac{\sqrt{2}}{2}$, for some positive integers $d_i$ (see Figure \ref{DouglasSchroder}).  The above set-up of drawn-in diagonal gives a new lattice whose fundamental regions are unit squares or triangles (halves of a unit square). We note that the triangles only appear along the drawn-in diagonals.

Given a positive integer $a$, we define the region $D_a(d_1,\dotsc,d_k)$ as follows (see Figure \ref{DouglasSchroder} for an example).  The southwestern and northeastern boundaries of $D_{a}(d_1,\dots,d_k)$ are defined in the next paragraph.

Color the new lattice black and white so that any two fundamental regions sharing an edge have opposite color, and  the fundamental regions passed through by $\ell$ are white.  Starting from a lattice point $A$ on $\ell$, we take unit steps south or east so that for each step the fundamental region on the right is black. We meet $\ell'$ at another lattice point $B$; and the described path from $A$ to $B$ is the northeastern boundary of our region. Next, we pick the lattice point $D$ on $\ell$ to the left of $A$ so that the distance between $A$ and $D$ is  $a\sqrt{2}$. The southwestern boundary of our region is obtained by reflecting the northeastern  one about the perpendicular bisector of segment $AD$, and reversing the directions of its steps (from south to north, and from east to west). Let $C$ be the intersection of the southwestern boundary and $\ell'$.

Finally, we connect $D$ to $A$ and $B$ to $C$ by two zigzag lattice paths. These two zigzag lattice paths are the northwestern and southeastern boundaries, and they complete the boundary of the region $D_a(d_1,\dotsc,d_k)$. We call the resulting region a \textit{generalized Douglas region}, and the four lattice points $A$, $B$, $C$ and $D$ the \textit{vertices} of the regions.

As mentioned in \cite{Tri}, to eliminate disconnected regions and  regions with no tiling,   we assume in addition that our generalized Douglas regions have disjoint southwestern and northeastern boundaries, and that the diagonal $\ell'$ passes through white unit squares.

%

Following the language in the prequel \cite{Tri} of the paper, we call the fundamental regions in a generalized Douglas region \textit{cells}.  There are two types of cells, square and triangular cells; and the latter come in two orientations: they may point towards $\ell'$ or away from $\ell'$. We call them \textit{down-pointing triangles} or \textit{up-pointing triangles}, respectively.

A \textit{(southwest-to-northeast) cell-line} consists of all the triangular cells of a given color with bases resting on a fixed lattice diagonal, or of all the square cells passed through by a fixed lattice diagonal. Define the \textit{width} $w$ of our region to be the number of white squares along the bottom of the region.  The generalized Douglas region in Figure \ref{DouglasSchroder} has width $w=8$. A  \textit{regular} cell is a  black square or a black up-pointing triangle. We note that if a cell in a cell-line is regular, so are all cells in the cell-line. This cell-line is called a \emph{regular cell-line}.

The author proved in \cite{Tri} the following tiling formula of a generalized Douglas region by extending Eu and Fu's lattice path method \cite{Eu}.
\begin{thm}[Theorem 4 in \cite{Tri}]\label{gendoug}
Assume that  $a,$ $d_1,$ $\dotsc,$ $d_k$ are positive integers, for which the generalized Douglas region $D_{a}(d_1,\dotsc,d_k)$ has the width  $w$, and the western and eastern vertices (i.e. the vertices $B$ and $D$) on the same horizontal line. Then
\begin{equation}\label{gendougeq}
\M(D_a(d_1,\dotsc,d_k))=2^{\mathcal{C}-w(w+1)/2},
\end{equation}
where $\mathcal{C}$ is the number of regular cells of the region.
\end{thm}

One readily sees that the Aztec diamond theorem and Douglas' theorem are two special cases of Theorem \ref{gendoug}, corresponding to the case when $k=1$ and $a=d_1=n$, and the case when $k=2n\geq 2$, $d_1=d_k=1$, $a=k$, and $d_2=d_3=\dotsc=d_{k-1}=2$, respectively.

The goal of this paper is to present a new proof for the Theorem \ref{gendoug} using Kuo's graphical condensation \cite{Kuo}.
\section{New proof of Theorem \ref{gendoug}}

 A \textit{perfect matching} of a graph $G$ is a collection of disjoint edges covering all vertices of $G$. The \textit{dual graph} of a region $R$ is the graph whose vertices are the fundamental regions in $R$, and whose edges connect precisely two fundamental regions sharing an edge. The tilings of the region  $R$ can be identified with the perfect matchings of its dual graph. In the view of this, we use the notation $\M(G)$ for the number of perfect matchings of a graph $G$.

Eric H. Kuo (re-)proved the Aztec diamond theorem by using a method called ``graphical condensation" \cite{Kuo}. The key of his proof is the following combinatorial interpretation of the Desnanot-Jacobi identity in linear algebra (see e.g. \cite{Mui}, pp. 136--149).

\begin{thm} [Kuo  \cite{Kuo}] \label{kuothm} Let $G=(V_1,V_2,E)$ be a planar bipartite graph in which $|V_1|=|V_2|$. Assume that $x,y,z$ and $t$ are four vertices appearing in a cyclic order on a face of $G$ so that $x,z\in V_1$ and $y,t\in V_2$. Then
\begin{equation}\label{kuoeq}
\M(G)\M(G-\{x,y,z,t\})
=\M(G-\{x,y\})\M(G-\{z,t\})+\M(G-\{t,x\})\M(G-\{y,z\}).
\end{equation}
\end{thm}

\medskip

The $k-1$ drawn-in diagonals divide the generalized Douglas region $D_{a}(d_1,\dotsc,d_k)$ into $k$ parts called \textit{layers}. The first layer is the part above the top drawn-in diagonal, the last layer is the part below the bottom drawn-in diagonal, and the $i$-th layer (for $1<i<k$) is the part between the $(i-1)$-th and the $i$-th drawn-in diagonals. We call $d_i$ the \emph{thickness} of  the $i$-th layer.
The sum of all $d_i$'s is called the \emph{thickness} of the region.

\medskip

\begin{proof}[Proof of Theorem \ref{gendoug}]

We prove (\ref{gendougeq}) by induction on the thickness $h=\sum_{i=1}^{k}d_i$ of the generalized Douglas region $\mathcal{D}:=D_{a}(d_1,d_2,\dotsc,d_k)$.

The base cases are the situations when $h=1$ or $k=1$.

\medskip

If $k=1$, then our region is exactly the Aztec diamond region of order $a$, and (\ref{gendougeq}) follows from the Aztec diamond theorem \cite{Elkies, Elkies2}. If $h=1$, then $k=1$ and $d_1=1$, and our region is the Aztec diamond of order 1. Again, (\ref{gendougeq}) follows from the Aztec diamond theorem.

\begin{figure}\centering
\includegraphics[width=6cm]{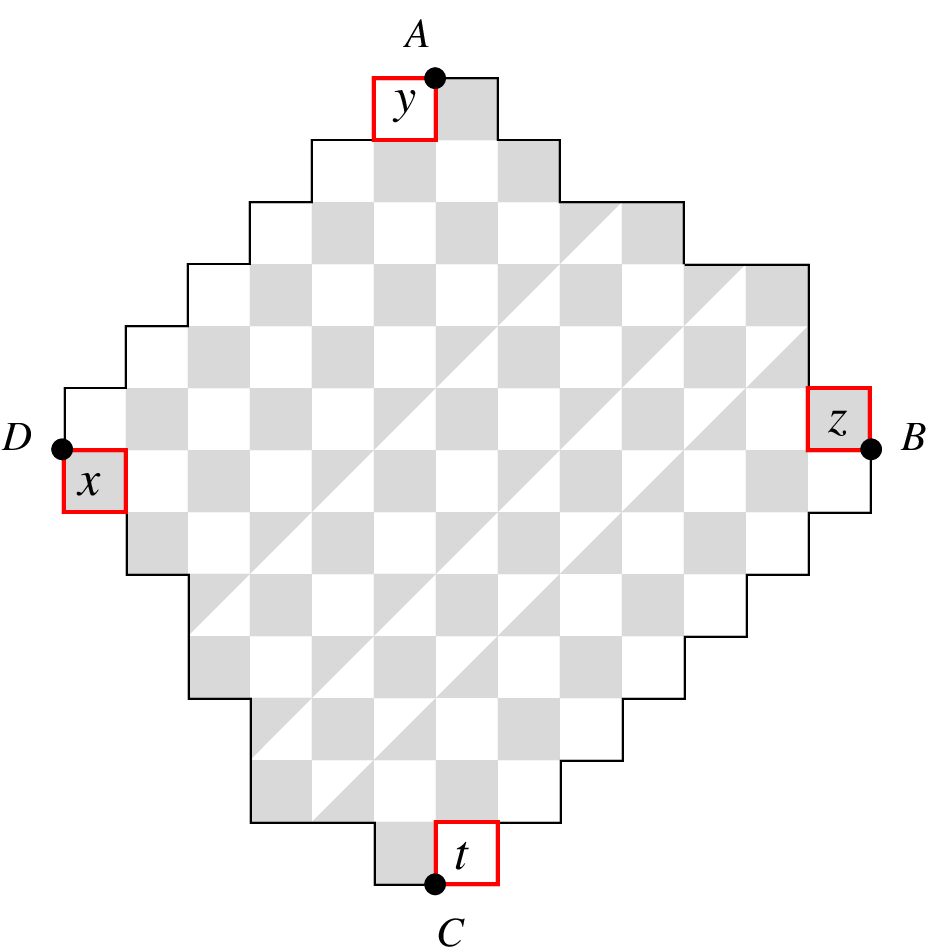}
\caption{How we apply Kuo condensation to the dual graph of a generalized Douglas region.}
\label{DouglasKuo}
\end{figure}

\medskip

Our induction step is based on Kuo's Condensation Theorem \ref{kuothm}.

For the induction step, we assume that $k,h>1$ and that  (\ref{gendougeq}) holds for any generalized Douglas regions in which the thickness is strictly less than $h$.

We apply Kuo's Condensation Theorem \ref{kuothm} to the dual graph $G$ of $\mathcal{D}$. Each vertex of $G$ corresponds to a cell of $\mathcal{D}$. The four vertices $x,y,z,t$ in the theorem are chosen as in Figure \ref{DouglasKuo}. In particular, the cells corresponding to $x$ and $z$ are respectively the black cells adjacent to the vertices $D$ and $B$; and $y$ and $t$ correspond to the white cells adjacent to the vertices $A$ and $C$, respectively.

\begin{figure}\centering
\includegraphics[width=13cm]{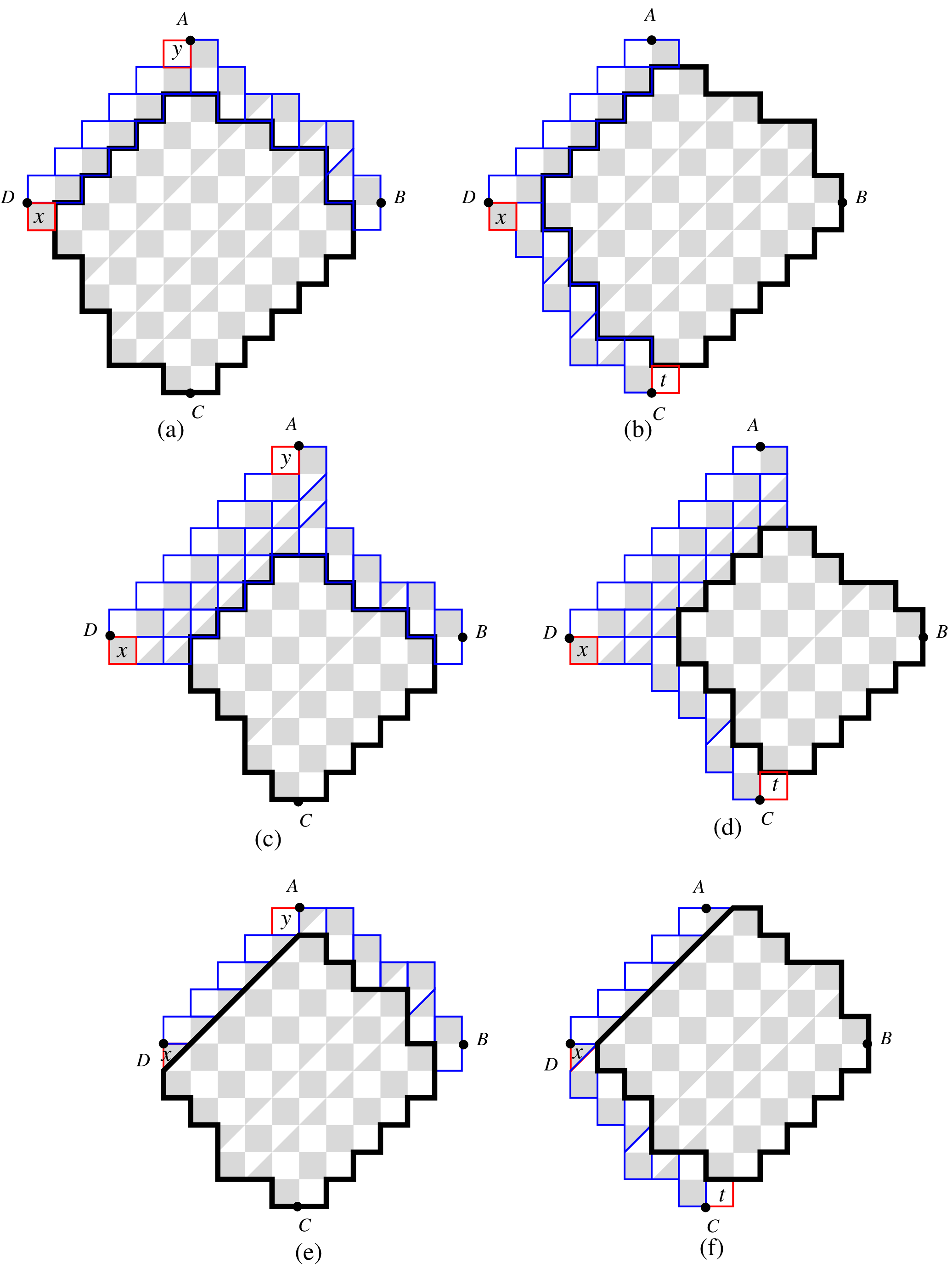}
\caption{Obtaining the equation $\M(G-\{x,y\})=\M(G-\{x,t\})=\M(\mathcal{D}_1)$ when $\mathcal{D}=D_{6}(5,3,2,3)$ ((a) and (b)), $\mathcal{D}=D_{7}(2,1,6,4)$ ((c) and (d)), and $\mathcal{D}=D_{5}(1,5,2,3)$ ((e) and (f)).}
\label{DouglasKuo3}
\end{figure}

Let us consider the region corresponding to the graph $G-\{x,y\}$ (see Figures \ref{DouglasKuo3}(a), (c) and (e)). There are several tiles that are forced to be in any tiling of the region. By removing these forced tiles, we get a new generalized Douglas region having the same number of tilings as the original one. In particular,   if $d_1\geq 3$, we have one stack of forced tiles running along the northwest side and another stack along the northeast side (see Figure \ref{DouglasKuo3}(a)). Removing these stacks, we get the region $D_{a-1}(d_1-2,d_2,\dotsc,d_k)$ (illustrated by the region restricted by the bold contour in Figure \ref{DouglasKuo3}(a)).
  If $d_2=2$, we may have more than one stack of forced tiles along the northwest side, together with a stack along the northeast side (illustrated by Figure \ref{DouglasKuo3}(c)).  In particular,  the removal of forced tiles gives us the region
$D_{a-m}(d_{m}-1,d_{m+1},\dotsc,d_k)$, where $m$ is the smallest index greater than 1 so that $d_m\geq 2$.
Finally, if $d_1=1$, we have one stack of forced tiles along each of the northeast and northwest sides.  The resulting region is obtained from the generalized Douglas region $D_{a}(d_{2},d_{3},\dotsc,d_k)$ by replacing the top row of white squares by a row of white triangles (shown in Figure \ref{DouglasKuo3}(e)).  However, this replacement does not change the number of tilings. More precisely, the dual graph of the resulting region is isomorphic to the dual graph of the region $D_{a}(d_{2},d_{3},\dotsc,d_k)$.
In summary, if we define a formal region $\mathcal{D}_1$ as
\begin{equation}
\mathcal{D}_1:=
\begin{cases}
D_{a-1}(d_1-2,d_2,\dotsc,d_k) &\text{if $d_1\geq 3$;}\\
D_{a-m}(d_{m}-1,d_{m+1},\dotsc,d_k) &\text{if $d_1=2$}\\
D_{a}(d_{2},d_{3},\dotsc,d_k) &\text{if $d_1=1$,}
\end{cases}
\end{equation}
then we always have
\begin{equation}\label{proof1eq1}
\M(G-\{x,y\})=\M(\mathcal{D}_1).
\end{equation}
Similarly, we also have
\begin{equation}\label{proof1eq2}
\M(G-\{t,x\})=\M(\mathcal{D}_1)
\end{equation}
(see Figures \ref{DouglasKuo3}(b), (d) and (f) for the cases $d_1\geq3$, $d_1=2$ and $d_1=1$, respectively.)

By symmetry, we have
\begin{equation}\label{proof1eq3}
\M(G-\{z,t\})=\M(G-\{y,z\})=\M(\mathcal{D}_2),
\end{equation}
where
\begin{equation}
\mathcal{D}_2:=
\begin{cases}
D_{a-1}(d_1,d_2,\dotsc,d_k-2) &\text{if $d_k\geq 3$;}\\
D_{a}(d_1,d_{2},\dotsc,d_q-1) &\text{if $d_k=2$;}\\
D_{a}(d_{2},d_{3},\dotsc,d_{k-1}) &\text{if $d_k=1$}
\end{cases}
\end{equation}
 and where $q$ is the largest index smaller than $k$ so that $d_q\geq 2$.

\medskip

Without loss of generality, we assume from now on that $d_1\leq d_k$ (otherwise, one can rotate our region $180^{\circ}$ to get a new generalized Douglas region with the bottom layer thicker than the top layer).

\medskip

There are three cases to distinguish, based on the value of $d_k$.

\bigskip

\begin{figure}\centering
\includegraphics[width=13cm]{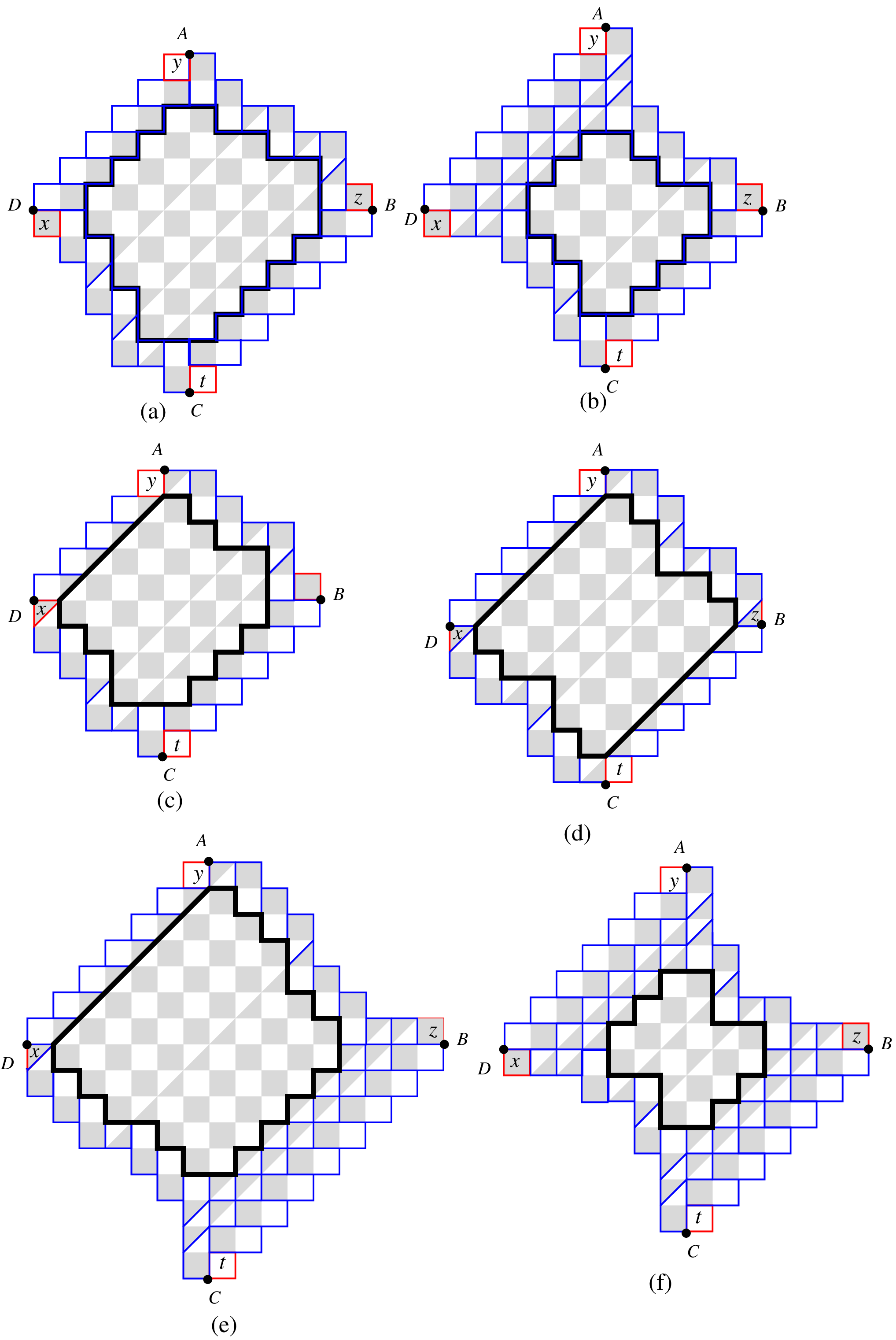}
\caption{Obtaining the equation $\M(G-\{x,y,z,t\})=\M(\mathcal{D}_3)$ when $\mathcal{D}=D_{6}(5,3,2,3)$ ((a) and (b)), $\mathcal{D}=D_{5}(1,5,2,3)$ ((c) and (d)), and $\mathcal{D}=D_{7}(1,6,6,1,2)$ ((e) and (f)).}
\label{DouglasKuo4}
\end{figure}

\textit{Case I.  $d_k\geq 3$.}

\medskip

We consider the region corresponding to the graph $G-\{x,y,z,t\}$. By removing all tiles forced by the four cells corresponding to $x,y,z,t$, we always get a new generalized Douglas region $\mathcal{D}_3$ defined by
\begin{equation}
\mathcal{D}_3:=
\begin{cases}
D_{a-1}(d_1-2,d_2,\dotsc,d_k-2) &\text{if $d_1\geq 3$;}\\
D_{a-m}(d_{m}-1,d_{m+1},\dotsc,d_k-2) &\text{if $d_1=2$;}\\
D_{a}(d_{2},d_{3},\dotsc,d_k-2) &\text{if $d_1=1$}\\
\end{cases}
\end{equation}
(see Figures \ref{DouglasKuo4}(a), (b) and (c) for the cases $d_1\geq 3$, $d_1=2$ and $d_1=1$).
Thus,  we get
\begin{equation}\label{proof1eq4}
\M(G-\{x,y,z,t\})=\M(\mathcal{D}_3).
\end{equation}
Plugging (\ref{proof1eq1})--(\ref{proof1eq4}) into the equation (\ref{kuoeq}) in Kuo's Theorem \ref{kuothm}, we get
\begin{equation}\label{proof1eq5}
\M(\mathcal{D})\M(\mathcal{D}_3)=2\M(\mathcal{D}_1)\M(\mathcal{D}_2).
\end{equation}
One readily sees that all the generalized Douglas regions $\mathcal{D}_i$ have thickness strictly less than $h$. By the induction hypothesis, their numbers of tilings are all given by (\ref{gendougeq}). Substituting these formulas into (\ref{proof1eq5}), we get
\begin{equation}\label{proof1eq6}
\M(\mathcal{D})=2^{\mathcal{C}_1+\mathcal{C}_2-\mathcal{C}_3-\binom{w_1+1}{2}-\binom{w_2+1}{2}+\binom{w_3+1}{2}+1},
\end{equation}
where $w_i$ and $\mathcal{C}_i$ are the width and the number of regular cells of $\mathcal{D}_i$, for $i=1,2,3$.

Comparing the widths of $\mathcal{D}$ and $\mathcal{D}_1$, we get
\begin{equation}\label{c1a}
w_1=w-1.
\end{equation}
Similarly, we obtain
\begin{equation}\label{c1b}
w_2=w-1\quad\text{ and  }\quad w_3=w-2.
\end{equation}
We now consider the regular cells of $\mathcal{D}$ which are not in $\mathcal{D}_1$. These regular cells only appear on the topmost black cell-line (which consists of  $a+1$ black squares) or at the left end of any other regular cell-lines. Since the total number of  regular cell-lines is equal to $w$ (see equation (4) in \cite[pp. 5]{Tri}), we obtain
\begin{equation}\label{c1c}
\mathcal{C}_1=\mathcal{C}-(a+1)-(w-1)=\mathcal{C}-a-w.
\end{equation}
Similarly, we have
\begin{equation}\label{c1d}
\mathcal{C}_2=\mathcal{C}-2w \quad
\text{ and }
 \quad \mathcal{C}_3=\mathcal{C}-a-3w+2.
\end{equation}
By (\ref{c1a})--(\ref{c1d}),  the equation  (\ref{gendougeq}) follows from (\ref{proof1eq6}).

\medskip

\quad\textit{Case II. $d_k=2$.}

\medskip

By the assumption $d_k\geq d_1$, we have $d_1$ equals $1$ or $2$. Applying the same process as in Case I, we obtain also
\begin{equation}
\M(\mathcal{D})\M(\mathcal{D}_3)=2\M(\mathcal{D}_1)\M(\mathcal{D}_2),
\end{equation}
where $\mathcal{D}_1$ and $\mathcal{D}_2$ are defined as in Case I, and
\begin{equation}
\mathcal{D}_3:=
\begin{cases}
D_{a}(d_{2},d_{3},\dotsc,d_{q}-1) &\text{if $d_1=1$;}\\
D_{a-m}(d_{m}-1,d_{m+1},\dotsc,d_{q}-1) &\text{if $d_1=2$}.
\end{cases}
\end{equation}
(See Figures \ref{DouglasKuo4}(e) and (f), respectively, for the cases $d_1=1$ and $d_1=2$.) All the above $\mathcal{D}_i$ regions have thickness less than $h$, so by the induction hypothesis, we still have the equation (\ref{proof1eq6}).

We also have  the two equalities (\ref{c1a}) and (\ref{c1c}) as in Case I. Moreover, one can verify the following facts:
\begin{equation}\label{c4a}
w_2=w-(k-q)-1 \quad \text{ and } \quad w_3=w-(k-q)-2
\end{equation}
(see Figure \ref{Douglascompare} for the first equality, and Figures \ref{DouglasKuo4}(e) and (f) for the second equality).

\begin{figure}\centering
\includegraphics[width=6cm]{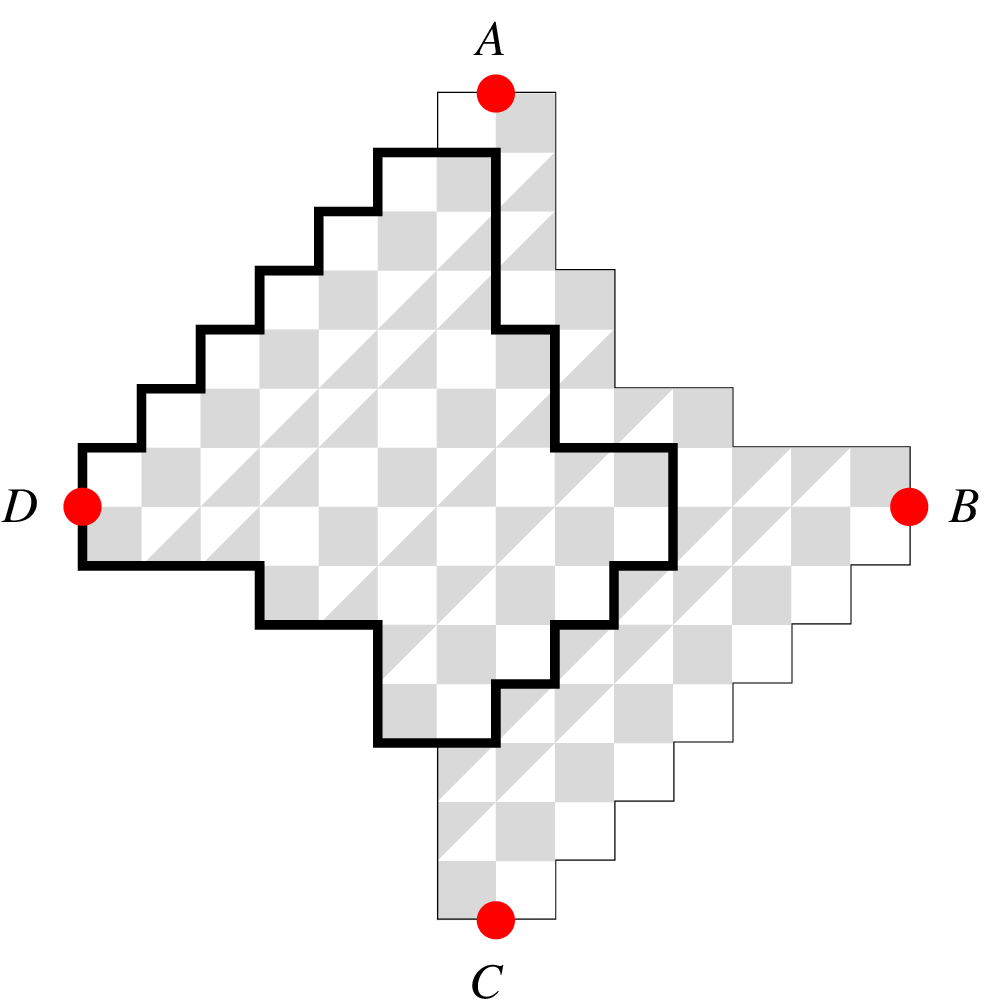}
\caption{Comparing the regions $\mathcal{D}=D_{7}(2,1,3,2,3,1,2)$ and $\mathcal{D}_2=D_{7}(2,1,3,2,2)$ (restricted by the bold contour).}
\label{Douglascompare}
\end{figure}

Next, let us consider the difference between the numbers of regular cells of $\mathcal{D}$ and $\mathcal{D}_2$, i.e. $\mathcal{C}-\mathcal{C}_2$. The difference counts all regular cells of $\mathcal{D}$, which are not in $\mathcal{D}_2$. These regular cells are in the last cell-line of black squares, in the $k-q$ cell-lines of black up-pointing triangles above it, or on the left end of each other regular cell-line  (see Figure \ref{Douglascompare}). The last cell-line of black squares contributes $w+1$  such regular cells, the $k-q$ cell-lines of up-pointing triangles above it contribute a total  of $\sum_{i=1}^{k-q}(w-i+1)$ regular cells, and the $w-(k-q+1)$ other regular cell-lines contribute 1 regular cell each. Thus, we have
\begin{equation}\label{c4b}
\mathcal{C}_2=\mathcal{C}-2w-\sum_{i=1}^{k-q}(w-i).
\end{equation}
Similarly, we get
\begin{equation}\label{c4c}
\mathcal{C}_3=\mathcal{C}-a-3w-(k-q)+2-\sum_{i=1}^{k-q}(w-i).
\end{equation}

Substituting (\ref{c1a}), (\ref{c1c}),  (\ref{c4a}), (\ref{c4b}) and (\ref{c4c}) into the equation (\ref{proof1eq6}) and working on simplifications, we  get (\ref{gendougeq}).
\medskip

\quad\textit{ Case III. $d_k=1$.}

\medskip

Since we are assuming $d_k\geq d_1$,  $d_1$ must be  $1$ in this case. Repeating  our machinery as  in the previous cases, we have also
\begin{equation}
\M(\mathcal{D})\M(\mathcal{D}_3)=2\M(\mathcal{D}_1)\M(\mathcal{D}_2),
\end{equation}
where $\mathcal{D}_1$ and  $\mathcal{D}_2$ are defined as in Case I,  and $\mathcal{D}_3:=D_{a-1}(d_{2},d_{3},\dotsc,d_{k-1})$. Again, by the induction hypothesis, we obtain the equation (\ref{proof1eq6}).

On the other hand, the equations (\ref{c1a}) and (\ref{c1c}) in Case I still hold; and it is not hard to see that
\begin{equation}\label{c5a}
w_2=w \quad\text{ and }\quad w_3=w-1,
\end{equation}
\begin{equation}\label{c5b}
\mathcal{C}_2=\mathcal{C}-w
\quad\text{ and } \quad
\mathcal{C}_3=\mathcal{C}-a-2w+1.
\end{equation}
Again, by plugging (\ref{c1a}), (\ref{c1c}), (\ref{c5a}) and (\ref{c5b}) into the equation (\ref{proof1eq6}),  we get (\ref{gendougeq}). This finishes our proof.
\end{proof}

\section{Concluding remarks}

This paper presents an example of the power of Kuo condensation, when applied to situations in which explicit conjectured formulas can be found. However, we may need to consider \emph{many} cases, especially, when our region has a complicated graphical structure as the generalized Douglas region.


Recall that the proof of Theorem \ref{gendoug} in this paper and  the previous proof in \cite{Tri} extend respectively the ideas of Kuo \cite{Kuo} and Fu and Eu \cite{Eu} in the case of the Aztec diamonds. This suggests that Theorem \ref{gendoug} may also be proven by generalizing the other proofs of the Aztec diamond theorem (listed in the introduction).

\end{document}